\documentclass{ajmart}

\issueinfo{2}{3}{September}{2008}

\copyrightinfo{2008}{Aulona Press  \emph{(Albanian J. Math.)}}

\def\issn{{\sc ISSN} 1930-1235: }
\def\issueyear{2008}

\PII{\issn  (\issueyear )}

\pagespan{199}{213}

\usepackage{amsmath,amssymb,amsbsy,amsfonts,amsthm,latexsym,amsopn,amstext, amsxtra,euscript,amscd, hyperref}
\usepackage{graphicx}
\usepackage{cite}

\usepackage{xy}
\input xy \xyoption{all}

\newtheorem{thm}{Theorem}[section]
\newtheorem{lemma}{Lemma}

\newtheorem{remark}{Remark}
\theoremstyle{definition}

\newcommand{\Aut}{\mbox {Aut}}
\newcommand{\ch}{\mbox {char }}

\newcommand{\bP}{{\mathbb P}^1}

\newcommand\M{\mathcal M}
\newcommand\C{\mathcal C}

\newcommand{\cH}{\mathcal H}

\newcommand\X{\mathcal X}
\newcommand\G{\bar G}
\def\si{\sigma}

\def\p{\phi}
\def\P{\Phi}
\def\d{\delta}
\def\U{\mathcal U}

\def\s{s}
\def\t{t}
\def\r{r}
\def\lar{\longrightarrow}

\begin{document}

\title[Equations of cyclic curves]{Determining equations of families of cyclic curves}

\thanks{Both authors were supported by a NATO grant, ICS. EAP. ASI No 982903}

\keywords{algebraic curves, Hurwitz spaces, equations}

\subjclass[2000]{Primary: 14Hxx, Secondary: 14H37, 14H10, }

\maketitle

\begin{center}
{\sc R. Sanjeewa } \\
\vspace*{1ex}
\emph{Department of Mathematics and Statistics\\
Oakland University, \\
Rochester, MI, 48309.\\
rsanjeew@oakland.edu\\}
\end{center}

\vspace*{3ex}

\begin{center}
{\sc T. Shaska} \\
\vspace*{1ex}
\emph{Department of Mathematics\\
University of Vlora, \\
Vlora, Albania\\
shaska@univlora.edu.al}
\end{center}

\begin{abstract}
In previous work we determined automorphism groups of cyclic algebraic curves defined over fields of any odd characteristic. In this paper we determine parametric equations of families of curves for each automorphism group for such curves.
\end{abstract}

\section{Introduction}
Let $\X_g$ be an algebraic curve of genus $g\geq2$ defined over a algebraically closed field of characteristic $p\neq2$. If an automorphism group of a algebraic curve has normal cyclic subgroup such that the quotient space has genus zero, then such a curve is called a \emph{cyclic curve}. We have studied automorphism groups of \emph{cyclic curves} in \cite{SS}, where we have listed all automorphism groups as well as ramification signatures of corresponding covers. In this paper we give a corresponding parametric equation for each family in \cite{SS}.

In the second section we briefly introduce basic facts on cyclic curves and their automorphism. Let $G=\Aut(\X_g)$ automorphism group of given cyclic curve $\X_g$, the reduced automorphism group is $\G:= \Aut(\X_g)/\langle w \rangle$, where $C_n=\left\langle w\right\rangle$ such that $g(\X^{C_n})=0$. This group $\G$ is embedded in $PGL_2(k)$ and therefore is isomorphic to one of $C_m $, $ D_m $, $ A_4 $, $ S_4 $, $ A_5 $, \emph{a semi direct product of elementary Abelian group  with  cyclic group}, $ PSL(2,q)$, or  $PGL(2,q)$. Then, $\G$ acts on a genus 0 field $k(x)$. We determine a rational function $\phi(x)$ that generates the fixed field $k(x)^{\G}$ in all cases cf. Lemma ~\ref{l1}.

In section three, we determine the ramification signature $\si$ of the cover $\P(x):\X_g \to \bP$ with monodromy group $G:=\Aut(\X_g)$. Moduli spaces of covers $\P$ are Hurwitz spaces, which we denoted by $\cH_{\si}$. There is a map $\P_{\si}:\cH_{\si} \to \M_g$, where $\M_g$ is the moduli space of genus $g$ algebraic curves. The image of this map is a subvariety of $\M_g$, which we denoted by $\cH(G,\si)$. The dimension of $\cH(G,\si)$ is determined.  Hence, we have 
\[ \X_g    \buildrel{C_n} \over \lar  \bP \buildrel {\G} \over \lar \bP \]
We list all possible automorphism groups, their signatures, and dimension of the loci $\cH(G,\si)$.

In the last section, we determine the equations of families of curves for a given group. Using the rational function $\phi(x)$ we are able to determine parametric equation of each family $\cH(G,\si)$. Since we know $\phi(x)$, we can find the branch points and then determine the equation of the curve from these branch points. We list corresponding equations of families of curves which we have listed in section three.

Throughout this paper we let  $g \geq$ 2 be a fixed integer, $\X$ a genus $g$ cyclic curve, $G=\Aut(\X)$ and $C_n\triangleleft G$ such that $g(\X^{C_n})=0$.

\section{Preliminaries}

Let $\X_g$ be genus $g\geq 2$ cyclic curve defined over an algebraically closed field $k$ of characteristic $p\neq 2$. We take the equation of $\X_g$ to be $y^n=F(x)$, where $deg(F)=2g+2$. Let $K:=k(x,y)$ be the function field of $\X_g$. Then $K$ is a degree $n$ extension field of $k(x)$ ramified exactly at $d=2g+2$ places $\alpha_1,...,\alpha_d$ of $k(x)$.

Let $G=Aut(K/k)$. Since $k(x)$ is the only genus 0 subfield of degree $n$ of $K$, then $G$ fixes $k(x)$. Thus $Gal(K/k(x))=\left\langle w\right\rangle$, with $w^n=1$. Then the group $\G:=G/\left\langle w\right\rangle$ is called \emph{reduced automorphism group}. By the theorem of Dickson, $\G$ is isomorphic to one of the following:  $C_m $, $ D_m $, $ A_4 $, $ S_4 $, $ A_5$, $ PSL(2,q)$ and $PGL(2,q)$,  or a  semi direct product  of elementary  Abelian  group  with  cyclic  group, defined as  
\[ K_m:=\left\langle\left\{ \si_a, \t |a \in \U_m\right\}\right\rangle, \, \textit{ where  }  \, \U_m :=\{a \in k| (a\prod_{j=0}^{\frac{p^t-1}{m}-1}(a^m-b_j))=0\} \]
and $\t(x)=\xi^2x, \quad \si_a(x)=x+a$, for each $a \in \U$, $b_j \in k^*$ and $\xi$ is a primitive $2m$-th root of unity; see \cite{VM}. $\U_m$ is a subgroup of the additive group of $k$.

The group $\G$ acts on $k(x)$ via the natural way. The fixed field is a genus 0 field, say $k(z)$. Thus $z$
is a degree $|\G|$ rational function in $x$, say $z=\p(x)$. The following lemma determines rational functions for all $\G$; see \cite{SS}.

Let $\phi_0:\X_g \to \bP$ and $\phi:\bP \to \bP$ be covers which correspond to the extensions $K/k(x)$ and $k(x)/k$ respectively. Then, $\psi:=\phi \circ \phi_0$ has monodromy group $G:=Aut(\X_g)$. By basic covering theory, the group $G$ is embedded in the group $S_l$, where $l=deg(\psi)$. There is an $r$-tuple $\overline{\si}:=(\si_1,...,\si_r)$, where $\si_i \in S_l$ such that $\si_1,...,\si_r$ generate $G$ and $\si_1...\si_r=1$. The signature of $\P$ is an $r$-tuple of conjugacy classes $\C:=(C_1,...,C_r)$ in $S_l$ such that $C_i$ is the conjugacy class of $\si_i$. We can find the signature of $\psi_0:\X_g \to \bP$ by using the signature of $\phi:\bP \to \bP$ and Riemann-Hurwitz formula.

Moduli spaces of covers $\psi$ are Hurwitz space, which we denoted by $\cH_{\si}$. There is a map $\P_{\si}:\cH_{\si} \to \M_g$, where $\M_g$ is the moduli space of genus $g$ algebraic curves. The image of this map is a subvariety of $\M_g$, which we denoted by $\cH(G,\si)$. Using the signature  of $\psi$ and Riemann-Hurwitz formula, one can find out dimension of $\cH(G,\si)$, which we denoted by $\d$.

We summarize all in the following Lemma: 

\begin{lemma}\label{l1}
Let $k$ be an algebraically closed field of characteristic $p$, $H_t$  a subgroup of the additive group of $k$ with $| H_t | = p^t$ and $b_j \in k^*$, and $\G$ be a finite subgroup of $PGL_2(k)$ acting on the field $k(x)$.
Then, $\G$ is isomorphic to one of the following groups $C_m$, $D_m$, $A_4$, $S_4$, $A_5$, $U= C_p^t$, $K_m$, $PSL_2(q)$ and $PGL_2(q)$, where $q=p^f$ and $(m,p)=1$. Moreover, the fixed subfield $k(x)^{\G}=k(z)$ is given by Table ~\ref{t1}, where $\alpha =\frac{q(q-1)}{2}, \quad \beta= \frac{q+1}{2}$.
\end{lemma}

\begin{table}[hb]
\begin{center}
\begin{tabular}{cccc}
$Case$ & $\G$ & $z$ & $Ramification$  \\
\hline \\
1 & $C_m$, $(m,p)=1$& $x^m$ & $(m,m)$\\ \\
2 & $D_{2m}$, $(m,p)=1$& $x^m+\frac{1}{x^m}$ & $(2,2,m)$\\ \\
3 & $A_4, \, p\neq 2, 3$ & $\frac{x^{12}-33x^8-33x^4+1}{x^2(x^4-1)^2}$ & $(2,3,3)$\\ \\
4 & $S_4, \, p\neq 2, 3$ & $\frac{(x^8+14x^4+1)^3}{108(x(x^4-1))^4}$ & $(2,3,4)$\\ \\
5 & $A_5, \,p\neq 2, 3, 5$ & $\frac{(-x^{20}+228x^{15}-494x^{10}-228x^5-1)^3}{(x(x^{10}+11x^5-1))^5}$ & $(2,3,5)$\\ \\
  & $A_5, \,p=3$ & $\frac{(x^{10}-1)^6}{(x(x^{10}+2ix^5+1))^5}$ & $(6,5)$\\ \\
6 & $U$ & $  \displaystyle{\prod_{a \in H_t}} (x+a)$ & $(p^t)$\\ \\
7 & $K_m$ & $(x   \displaystyle{\prod_{j=0}^{\frac{p^t-1}{m}-1} } (x^m-b_j))^m$ & $(mp^t,m)$ \\ \\
8 & $PSL(2,q), \,p\neq 2$ & $\frac{((x^q-x)^{q-1}+1)^{\frac{q+1}{2}}}{(x^q-x)^{\frac{q(q-1)}{2}}}$ & $(\alpha,\beta)$ \\ \\
9 & $PGL(2,q)$ & $\frac{((x^q-x)^{q-1}+1)^{q+1}}{(x^q-x)^{q(q-1)}}$ & $(2\alpha,2\beta)$ \\ \\
\end{tabular}
\caption{Rational functions correspond to each $\G$} \label{t1}
\end{center}
\end{table}

\section{Automorphism groups and their signatures of cyclic curves}
As above $\G:=G/G_0$, where $G_0:=Gal(k(x,y)/k(x))$. The following theorem determines  ramification signatures and dimensions of $\d$ of $\cH(G,\si)$ for all $\G$ when $p>5$; see \cite{SS} for details.
\begin{thm} \label{th1}
The signature of cover $\P(x):\X \to \X^G$ and dimension $\delta$ is given in Table~\ref{t2}. In Table ~\ref{t2}, $m=|PSL_2(q)|$ for cases 38-41 and $m=|PGL_2(q)|$ for cases 42-45.
\end{thm}

\vspace{6ex}

\tiny

\begin{table}[h]
\begin{center}
\begin{tabular}{|c|c|c|c|c|}
\hline \hline
$\#$ & $\G$ & $\delta(G,C)$ & $\delta,n,g$ & $\C=(C_1,...,C_r)$ \\
\hline \hline
$1$ &$(p,m)=1$ & $\frac{2(g+n-1)}{m(n-1)}-1$ & $n<g+1$ & $(m,m,n,...,n)$  \\
$2$ & $C_m$ & $\frac{2g+n-1}{m(n-1)}-1$ & & $(m,mn,n,...,n)$ \\
$3$ &  & $\frac{2g}{m(n-1)}-1$ & $n<g$ &  $(mn,mn,n,...,n)$ \\
\hline \hline
$4$ & $(p,m)=1$ & $\frac{g+n-1}{m(n-1)}$ & & $(2,2,m,n,...,n)$  \\
$5$ &  & $\frac{2g+m+2n-nm-2}{2m(n-1)}$ & & $(2n,2,m,n,...,n)$  \\
$6$ & $D_{2m}$ & $\frac{g}{m(n-1)}$ & & $(2,2,mn,n,...,n)$ \\
$7$ & &$\frac{g+m+n-mn-1}{m(n-1)}$ & $n<g+1$ & $(2n,2n,m,n,...,n)$  \\
$8$ &  & $\frac{2g+m-mn}{2m(n-1)}$ & $g\neq 2$ & $(2n,2,mn,n,...,n)$  \\
$9$ &  & $\frac{g+m-mn}{m(n-1)}$ & $n<g$ & $(2n,2n,mn,n,...,n)$  \\
\hline \hline
$10$ &  & $\frac{n+g-1}{6(n-1)}$ & & $(2,3,3,n,...,n)$  \\
$11$ &$ A_4$ & $\frac{g-n+1}{6(n-1)}$ & & $(2,3n,3,n,...,n)$ \\
$12$ &  & $\frac{g-3n+3}{6(n-1)}$ & & $(2,3n,3n,n,...,n)$  \\
$13$ &  & $\frac{g-2n+2}{6(n-1)}$ & $\delta \neq 0$ & $(2n,3,3,n,...,n)$  \\
$14$ &  & $\frac{g-4n+4}{6(n-1)}$ & & $(2n,3n,3,n,...,n)$  \\
$15$ &  & $\frac{g-6n+6}{6(n-1)}$ & $\delta \neq 0$ & $(2n,3n,3n,n,...,n)$  \\
\hline \hline
$16$ &  & $\frac{g+n-1}{12(n-1)}$ & & $(2,3,4,n,...,n)$  \\
$17$ &  & $\frac{g-3n+3}{12(n-1)}$ & & $(2,3n,4,n,...,n)$  \\
$18$ &  &  $\frac{g-2n+2}{12(n-1)}$ & & $(2,3,4n,n,...,n)$ \\
$19$ &  & $\frac{g-6n+6}{12(n-1)}$ & & $(2,3n,4n,n,...,n)$  \\
$20$ &  $S_4$ & $\frac{g-5n+5}{12(n-1)}$ & & $(2n,3,4,n,...,n)$  \\
$21$ &  & $\frac{g-9n+9}{12(n-1)}$ & & $(2n,3n,4,n,...,n)$  \\
$22$ &  & $\frac{g-8n+8}{12(n-1)}$ & & $(2n,3,4n,n,...,n)$  \\
$23$ &  & $\frac{g-12n+12}{12(n-1)}$ & & $(2n,3n,4n,n,...,n)$  \\
\hline \hline
$24$ &  & $\frac{g+n-1}{30(n-1)}$ & & $(2,3,5,n,...,n)$  \\
$25$ &  & $\frac{g-5n+5}{30(n-1)}$ & & $(2,3,5n,n,...,n)$ \\
$26$ &  & $\frac{g-15n+15}{30(n-1)}$ & & $(2,3n,5n,n,...,n)$  \\
$27$ &  & $\frac{g-9n+9}{30(n-1)}$ & & $(2,3n,5,n,...,n)$  \\
$28$ &  $A_5$ & $\frac{g-14n+14}{30(n-1)}$ & & $(2n,3,5,n,...,n)$  \\
$29$ &  & $\frac{g-20n+20}{30(n-1)}$ & & $(2n,3,5n,n,...,n)$  \\
$30$ &  & $\frac{g-24n+24}{30(n-1)}$ & & $(2n,3n,5,n,...,n)$  \\
$31$ &  & $\frac{g-30n+30}{30(n-1)}$ & & $(2n,3n,5n,n,...,n)$  \\
\hline \hline
$32$ &  & $\frac{2g+2n-2}{p^t(n-1)}-2$ & & $(p^t,n,...,n)$  \\
$33$ & $U$ & $\frac{2g+np^{t}-p^t}{p^t(n-1)}-2$ & $(n,p)=1,n|p^t-1$ & $(np^t,n,...,n)$ \\
\hline \hline
$34$ &  & $\frac{2(g+n-1)}{mp^t(n-1)}-1$ & $(m,p)=1,m|p^t-1$ & $(mp^t,m,n,...,n)$  \\
$35$ &  & $\frac{2g+2n+p^t-np^t-2}{mp^t(n-1)}-1$ & $(m,p)=1,m|p^t-1$ & $(mp^t,nm,n,...,n)$  \\
$36$ & $K_m$ & $\frac{2g+np^t-p^{t}}{mp^t(n-1)}-1$ & $(nm,p)=1,nm|p^t-1$ & $(nmp^t,m,n,...,n)$\\
$37$ &  & $\frac{2g}{mp^t(n-1)}-1$ & $(nm,p)=1,nm|p^t-1$ & $(nmp^t,nm,n,...,n)$ \\
\hline \hline
$38$& & $\frac{2(g+n-1)}{m(n-1)}-1$ & $\left(\frac{q-1}{2},p\right)=1$ & $(\alpha,\beta,n,...,n)$  \\
$39$& $PSL_2(q)$  & $\frac{2g+q(q-1)-n(q+1)(q-2)-2}{m(n-1)}-1$ & $\left(\frac{q-1}{2},p\right)=1$ & $(\alpha,n\beta,n,...,n)$  \\
$40$& & $\frac{2g+nq(q-1)+q-q^2}{m(n-1)}-1$ & $\left(\frac{n(q-1)}{2},p\right)=1$ & $(n\alpha,\beta,n,...,n)$  \\
$41$& & $\frac{2g}{m(n-1)}-1$ & $\left(\frac{n(q-1)}{2},p\right)=1$ & $(n\alpha,n\beta,n,...,n)$  \\
\hline \hline
$42$& & $\frac{2(g+n-1)}{m(n-1)}-1$ & $(q-1,p)=1$ & $(2\alpha,2\beta,n,...,n)$  \\
$43$& $PGL_2(q)$  & $\frac{2g+q(q-1)-n(q+1)(q-2)-2}{m(n-1)}-1$ & $(q-1,p)=1$ & $(2\alpha,2n\beta,n,...,n)$  \\
$44$& & $\frac{2g+nq(q-1)+q-q^2}{m(n-1)}-1$ & $(n(p-1),p)=1$ & $(2n\alpha,2\beta,n,...,n)$ \\
$45$& & $\frac{2g}{m(n-1)}-1$ & $(n(q-1),p)=1$ & $(2n\alpha,2n\beta,n,...,n)$ \\ \\
\hline \hline
\end{tabular}
\medskip
\caption{The signature of curves and dimensions $\delta$ for $\ch >5 $} \label{t2}
\end{center}
\end{table}

\normalsize

\newpage

\begin{remark} \label{r1}
The above theorem gives signatures and dimensions for $p>5$. We know that $\G \cong C_m, D_m, A_4, S_4, U, K_m, PSL(2,q), PGL(2,q)$ when $p=5$ and $\G \cong C_m, D_m, A_5, U, K_m, PSL(2,q),PGL(2,q)$ when $p=3$; see \cite{VM}. All cases except $\G \cong A_5$ have ramification as $p>5$. Hence signatures and dimensions are the same as $p>5$. However, $\G \cong A_5$ has different ramification. Hence, that case has signatures and dimensions as in Table~\ref{t3}.
\end{remark}

\begin{table}[h]
\begin{center}
\begin{tabular}{|c|c|c|c|}
\hline \hline
$Case$ & $\G$ & $\delta(G,C)$ & $\C=(C_1,...,C_r)$ \\
\hline \hline
$a$ &  & $\frac{g+n-1}{30(n-1)}-1$ & $(6,5,n,...,n)$  \\
$b$ &  & $\frac{g+5n-5}{30(n-1)}-1$ & $(6,5n,n,...,n)$  \\
$c$ & $A_5$  & $\frac{g+6n-6}{30(n-1)}-1$ & $(6n,5,n,...,n)$ \\
$d$ & & $\frac{g}{30(n-1)}-1$ & $(6n,5n,n,...,n)$  \\
\hline \hline
\end{tabular}
\medskip
\caption{The signature of curve and dimension $\delta$ for $\G \cong A_5$, $p=3$}\label{t3}
\end{center}
\end{table}

The following theorem determines the list of all automorphism groups of cyclic algebraic curves defined over any
algebraically closed field of characteristic $p \neq 2$, details will be provided in \cite{SS}.

\begin{thm}\label{th2}
Let $\X_g$ be a genus $g\geq2$ irreducible cyclic curve defined over an algebraically closed field $k$
of characteristic $\ch (k) =p$, $G=Aut(\X_g)$, and $\G$ its reduced automorphism group. If $|G|>1 $ then is
$G$ is one of the following:

\begin{enumerate}

\item $\G \cong C_m$: Then, $G \cong C_{mn}$ or $\left\langle \r, \s \right|\r^n=1,\s^m=1,\s\r\s^{-1}=\r^l
\rangle$,  $(l,n)=1$ and $l^m\equiv 1$ (mod n).\\

\item If $\G \cong D_{2m}$ then $G \cong D_{2m} \times C_n$ or
\begin{align*}
\begin{split}
G'_{4}=&\left\langle \r, \s, \t \right|\r^n=1,\s^2=1,\t^2=1,(\s\t)^m=1,\s\r\s^{-1}=\r^l,\t\r\t^{-1}=\r^l \rangle \\
G'_{7}=& \left\langle \r, \s, \t \right|\r^n=1,\s^2=\r^{\frac{n}{2}},\t^2=\r^{\frac{n}{2}},(\s\t)^m=1,\s\r\s^{-1}=\r^l,\t\r\t^{-1}=\r^l \rangle \\
\end{split}
\end{align*}
where $(l,n)=1$ and $l^2\equiv 1$ (mod n) or
\begin{align*}
\begin{split}
G_{4}=& \left\langle \r, \s, \t \right|\r^n=1,\s^2=1,\t^2=1,(\s\t)^m=1,\s\r\s^{-1}=\r^l,\t\r\t^{-1}=\r^k \rangle \\
G_{5}=& \left\langle \r, \s, \t \right|\r^n=1,\s^2=\r^{\frac{n}{2}},\t^2=1,(\s\t)^m=1,\s\r\s^{-1}=\r^l,\t\r\t^{-1}=\r^k \rangle \\
G_{6}=& \left\langle \r, \s, \t \right|\r^n=1,\s^2=1,\t^2=1,(\s\t)^m=\r^{\frac{n}{2}},\s\r\s^{-1}=\r^l,\t\r\t^{-1}=\r^k \rangle \\
G_{7}=& \left\langle \r, \s, \t \right|\r^n=1,\s^2=\r^{\frac{n}{2}},\t^2=\r^{\frac{n}{2}},(\s\t)^m=1,\s\r\s^{-1}=\r^l,\t\r\t^{-1}=\r^k \rangle \\
G_{8}=& \left\langle \r, \s, \t \right|\r^n=1,\s^2=\r^{\frac{n}{2}},\t^2=1,(\s\t)^m=\r^{\frac{n}{2}},\s\r\s^{-1}=\r^l,\t\r\t^{-1}=\r^k \rangle \\
G_{9}=& \left\langle \r, \s, \t
\right|\r^n=1,\s^2=\r^{\frac{n}{2}},\t^2=\r^{\frac{n}{2}},(\s\t)^m=\r^{\frac{n}{2}},\s\r\s^{-1}=\r^l,\t\r\t^{-1}=\r^k
\rangle
\end{split}
\end{align*}
where $(l,n)=1$ and $l^2\equiv 1$ (mod n), $(k,n)=1$ and $k^2 \equiv 1$ (mod n).\\

\item If $\G \cong A_4$ and $p\neq2,3$ then $G \cong A_4 \times C_n$ or
\begin{align*}
\begin{split}
G'_{10}=& \left\langle \r, \s, \t \right|\r^n=1,\s^2=1,\t^3=1,(\s\t)^3=1,\s\r\s^{-1}=\r,\t\r\t^{-1}=\r^l \rangle\\
G'_{12}=& \left\langle \r, \s, \t \right|\r^n=1,\s^2=1,\t^3=\r^{\frac{n}{3}},(\s\t)^3=\r^{\frac{n}{3}},\s\r\s^{-1}=\r,\t\r\t^{-1}=\r^l \rangle\\
\end{split}
\end{align*}
where $(l,n)=1$ and $l^3\equiv 1$ (mod n) or
\begin{center}
$\left\langle \r, \s, \t
\right|\r^n=1,\s^2=\r^{\frac{n}{2}},\t^3=\r^{\frac{n}{2}},(\s\t)^5=\r^{\frac{n}{2}},\s\r\s^{-1}=\r,\t\r\t^{-1}=\r
\rangle , $ or
\end{center}
\begin{align*}
\begin{split}
G_{10}=& \left\langle \r, \s, \t \right|\r^n=1,\s^2=1,\t^3=1,(\s\t)^3=1,\s\r\s^{-1}=\r,\t\r\t^{-1}=\r^k \rangle\\
G_{13}=& \left\langle \r, \s, \t \right|\r^n=1,\s^2=\r^{\frac{n}{2}},\t^3=1,(\s\t)^3=1,\s\r\s^{-1}=\r,\t\r\t^{-1}=\r^k \rangle\\
\end{split}
\end{align*}
where $(k,n)=1$ and $k^3\equiv 1$ (mod n).\\

\item If $\G \cong S_4$ and $p\neq2,3$ then $G \cong S_4 \times C_n$ or
\begin{align*}
\begin{split}
G_{16}=& \left\langle \r, \s, \t \right|\r^n=1,\s^2=1,\t^3=1,(\s\t)^4=1,\s\r\s^{-1}=\r^l,\t\r\t^{-1}=\r \rangle\\
G_{18}=& \left\langle \r, \s, \t \right|\r^n=1,\s^2=1,\t^3=1,(\s\t)^4=\r^{\frac{n}{2}},\s\r\s^{-1}=\r^l,\t\r\t^{-1}=\r \rangle\\
G_{20}=& \left\langle \r, \s, \t \right|\r^n=1,\s^2=\r^{\frac{n}{2}},\t^3=1,(\s\t)^4=1,\s\r\s^{-1}=\r^l,\t\r\t^{-1}=\r \rangle \\
G_{22}=& \left\langle \r, \s, \t \right|\r^n=1,\s^2=\r^{\frac{n}{2}},\t^3=1,(\s\t)^4=\r^{\frac{n}{2}},\s\r\s^{-1}=\r^l,\t\r\t^{-1}=\r \rangle\\
\end{split}
\end{align*}
where $(l,n)=1$ and $l^2\equiv 1$ (mod n).\\

\item If $\G \cong A_5$ and $p\neq2,5$ then $G \cong A_{5}\times C_{n}$ or
\[\left\langle \r, \s, \t
\right|\r^n=1,\s^2=\r^{\frac{n}{2}},\t^3=\r^{\frac{n}{2}},(\s\t)^5=\r^{\frac{n}{2}},\s\r\s^{-1}=\r,\t\r\t^{-1}=\r
\rangle \]

\item If $\G \cong U$ then $G \cong U \times C_n$ or
\begin{equation*}
\left\langle \r,\s_1,\s_2,...,\s_t|\r^n=\s_1^p=\s_2^p=...=\s_t^p=1, \s_i\s_j=\s_j\s_i,
\s_i\r\s_i^{-1}=\r^{l}, 1\leq i,j\leq t\right\rangle
\end{equation*}
where $(l,n)=1$ and $l^p \equiv 1$ (mod n).\\

\item If $\G \cong K_m$ then $G \cong < \r,\s_1,...,\s_t,v|\r^n=\s_1^p=...=\s_t^p=v^m=1, \s_i\s_j=\s_j\s_i,
v\r v^{-1}=\r, \s_i\r\s_i^{-1}=\r^{l}, \s_iv\s_i^{-1}=v^{k}, 1\leq i,j\leq t >$  where $(l,n)=1$ and $l^p
\equiv 1$ (mod n), $(k,m)=1$ and $k^p \equiv 1$ (mod m) or
\begin{align*}
\begin{split}
G_{35}=&\left\langle \r,\s_1,...,\s_t|\r^{nm}=\s_1^p=...=\s_t^p=1, \s_i\s_j=\s_j\s_i,
\s_i\r\s_i^{-1}=\r^{l}, 1\leq i,j\leq t\right\rangle
\end{split}
\end{align*}
where $(l,nm)=1$ and $l^p \equiv 1 (\mod nm)$.\\

\item If $\G \cong PSL_{2}(q)$ then $G\cong PSL_2(q)\times C_n$ or $SL_2(3)$.\\

\item If $\G \cong PGL(2,q)$ then $G \cong PGL(2,q) \times C_n$.\\
\end{enumerate}
\end{thm}

\proof See \cite{SS}. \qed

\section{Equations of curves}

The group $\G$ is the monodromy group of the cover $\phi:\bP \to \bP$ with signature $(\si_1,\si_2,\si_3)$ as in section 2. We fix coordinates in $\bP$ as $x$ and $z$ respectively and from now on we denote the cover $\phi:\bP_x\to \bP_z$. Thus, $z$ is a rational function in $x$ of the degree $|\G|$. We denote by $q_1,q_2,q_3$ corresponding branch points of $\phi$. Let $S$ be the set of branch points of $\Phi:\X_g\to \bP_z$. Clearly $q_1,q_2,q_3 \in S$. Let $y^n=f(x)$ be the equation of $\X_g$ and $W$ be the images in $\bP_x$ of roots of $f(x)$ and \[V:=\bigcup_{i=1}^3 \phi^{-1}(q_i).\]

Let \[z=\frac{\Psi(x)}{\Upsilon(x)}, \textit{ where }     \Psi(x), \Upsilon(x) \in k[x].\]
Then we have \[z-q_i=\frac{\Gamma(x)}{\Upsilon(x)}\] for each branch point $q_i$, $i=1,2,3$, where $\Gamma(x) \in k[x]$. Hence,
\[\Gamma(x)=\Psi(x)-q_i\cdot\Upsilon(x)\] 
is degree $|\G|$ equation and multiplicity of all roots of $\Gamma(x)$ correspond to the ramification index for each $q_i$. Now we define the following three functions:

\begin{align} \label{e3}
\begin{split}
\varphi^r(x):=\Psi(x)-q_1 \cdot \Upsilon(x)\\
\chi^s(x):=\Psi(x)-q_2 \cdot \Upsilon(x)\\
\psi^t(x):=\Psi(x)-q_3 \cdot \Upsilon(x)
\end{split}
\end{align}
where superscript denote the ramification index of $q_i$.
Clearly, $\p^{-1}(S \backslash \{q_1,q_2,q_3\})\subset W $. Let $\lambda \in $ S  $\backslash$  $\{q_1,q_2,q_2\}$. The points in the fiber $\phi^{-1}(\lambda)$ are the roots of the equation:
\begin{equation} \label{e4}
\begin{split}
\Psi(x)-\lambda \cdot \Upsilon(x)=0
\end{split}
\end{equation}

Let
\begin{equation}\label{e5}
\begin{split}
G(x):= \prod_{\lambda \in  S  \backslash  \{q_1,q_2,q_2\}}(\Psi(x)-\lambda \cdot \Upsilon(x))
\end{split}
\end{equation}

There are following cases and corresponding equations of the curve $y^n=f(x)$ for each fixed $\phi$.

\medskip

\begin{center}
\begin{tabular}{lll}
& Intersection & $f(x)$\\
\hline 
&& \\
1) &$V \cap W = \emptyset $ & $G(x)$\\
2) &$V \cap W = \phi^{-1}(q_1) $ & $\varphi(x)\cdot G(x)$\\
3) &$V \cap W = \phi^{-1}(q_2) $ & $\chi(x)\cdot G(x)$\\
4) &$V \cap W = \phi^{-1}(q_3) $ & $\psi(x)\cdot G(x)$\\
5) &$V \cap W = \phi^{-1}(q_1)\cup \phi^{-1}(q_2) $ & $\varphi(x)\cdot \chi(x)\cdot G(x)$\\
6) &$V \cap W = \phi^{-1}(q_2)\cup \phi^{-1}(q_3) $ & $\chi(x)\cdot \psi(x)\cdot G(x)$\\
7) &$V \cap W = \phi^{-1}(q_1)\cup \phi^{-1}(q_3) $ & $\varphi(x)\cdot \psi(x)\cdot G(x)$\\
8) &$V \cap W = \phi^{-1}(q_1)\cup \phi^{-1}(q_2)\cup \phi^{-1}(q_3) $ & $\varphi(x)\cdot \chi(x)\cdot \psi(x)\cdot G(x)$\\
\end{tabular}
\end{center}

\medskip 

The following theorem gives us equations of families of curves for automorphism groups which are related to Theorem ~\ref{th1} and Theorem ~\ref{th2}.

\begin{thm}\label{th3}
Let $\X_g$ be a genus $g \geq$ 2 cyclic curve with $Aut(\X_g)=G$, where $G$ is related to the cases 1-45 in Table ~\ref{t2}. Then $\X_g$ has an equation as cases 1-45 in Table ~\ref{t6}.
\end{thm}

\begin{table}
\begin{center}
\begin{tabular}{|c|c|c|}
\hline
$\#$ &$\G$& $y^n=f(x)$ \\
\hline
1 &     & $x^{m\d}+a_1x^{m(\d-1)}+ \dots +a_\d x^{m}+1$\\
2 &$C_m$& $x^{m\d}+a_1x^{m(\d-1)}+ \dots +a_\d x^{m}+1$\\
3 &     & $x(x^{m\d}+a_1x^{m(\d-1)}+ \dots +a_\d x^{m}+1)$\\
\hline
4 &        & $F(x):= \prod_{i=1}^\d(x^{2m}+\lambda_ix^m+1)$\\
5 &        & $(x^m-1)\cdot F(x)$\\
6 &        & $x\cdot F(x)$\\
7 &$D_{2m}$& $(x^{2m}-1)\cdot F(x)$\\
8 &        & $x(x^m-1)\cdot F(x)$\\
9 &        & $x(x^{2m}-1)\cdot F(x)$\\
\hline
10 &     & $G(x):= \prod_{i=1}^\delta(x^{12}-\lambda_ix^{10}-33x^8+2\lambda_ix^6-33x^4-\lambda_ix^2+1)$\\
11 &     & $(x^4+2i\sqrt{3}x^2+1)\cdot G(x)$\\
12 &$A_4$& $(x^8+14x^4+1)\cdot G(x)$\\
13 &     & $x(x^4-1)\cdot G(x)$\\
14 &     & $x(x^4-1)(x^4+2i\sqrt{3}x^2+1)\cdot G(x)$\\
15 &     & $x(x^4-1)(x^8+14x^4+1)\cdot G(x)$\\
\hline
16 &     & $M(x)$\\
17 &     & $S(x)\cdot M(x)$\\
18 &     & $T(x)\cdot M(x)$\\
19 &     & $S(x)\cdot T(x)\cdot M(x)$\\
20 &$S_4$& $R(x)\cdot M(x)$\\
21 &     & $R(x)\cdot S(x)\cdot M(x)$\\
22 &     & $R(x)\cdot T(x)\cdot M(x)$\\
23 &     & $R(x)\cdot S(x)\cdot T(x).M(x)$\\
\hline
24 &     & $\Lambda(x)$\\
25 &     & $(x(x^{10}+11x^5-1))\cdot \Lambda(x)$\\
26 &     & $(x^{20}-228x^{15}+494x^{10}+228x^5+1)(x(x^{10}+11x^5-1))\cdot \Lambda(x)$\\
27 &     & $(x^{20}-228x^{15}+494x^{10}+228x^5+1)\cdot \Lambda(x)$\\
28 &$A_5$& $ Q (x)\cdot \Lambda(x)$\\
29 &     & $x(x^{10}+11x^5-1).\psi(x)\cdot \Lambda(x)$\\
30 &     & $(x^{20}-228x^{15}+494x^{10}+228x^5+1)\cdot \psi(x)\cdot\Lambda(x)$\\
31 &     & $(x^{20}-228x^{15}+494x^{10}+228x^5+1)(x(x^{10}+11x^5-1))\cdot\psi(x)\cdot\Lambda(x)$\\
\hline
32 &$U$& $B(x)$\\
33 &   & $B(x)$\\
\hline
34 &     & $\Theta(x)$\\
35 &$K_m$& $x\prod_{j=1}^{\frac{p^t-1}{m}}\left(x^m-b_j\right)\cdot\Theta(x)$\\
36 &     & $\Theta(x)$\\
37 &     & $x\prod_{j=1}^{\frac{p^t-1}{m}}\left(x^m-b_j\right)\cdot\Theta(x)$\\
\hline
38 &          & $\Delta(x)$\\
39 &$PSL_2(q)$& $((x^q-x)^{q-1}+1)\cdot\Delta(x)$\\
40 &          & $(x^q-x)\cdot\Delta(x)$\\
41 &          & $(x^q-x)((x^q-x)^{q-1}+1)\cdot\Delta(x)$\\
\hline
42 &          & $\Omega(x)$\\
43 &$PGL_2(q)$& $((x^q-x)^{q-1}+1)\cdot\Omega(x)$\\
44 &          & $(x^q-x)\cdot\Omega(x)$\\
45 &          & $(x^q-x)((x^q-x)^{q-1}+1)\cdot\Omega(x)$\\
\hline
\end{tabular}
\medskip
\caption{The equations of the curves related to the cases in Table ~\ref{t2}}\label{t6}
\end{center}
\end{table}

\noindent \textbf{Proof: }  We consider all cases one by one for the reduced automorphism group $\G$. 

\subsection{$\G \cong C_m$} Then, $\phi:\bP \to \bP$ has signature $(m,m)$. We identify the branch points of $\phi$ are $0$ and $\infty$. Let $q_1=\infty$, $q_2=0$. By Lemma ~\ref{l1}, we know that $\phi(x)=x^m$. Hence $\varphi(x)=1$ and $\chi(x)=x$. Let $\lambda_i \in $ S  $\backslash$  $\{0,\infty\}$. The points in the fiber $\phi^{-1}(\lambda_i)$ are the roots of the polynomial
\begin{equation*}
\begin{split}
G_{\lambda_i}(x):= x^m-\lambda_i
\end{split}
\end{equation*}
Now we can compute equations for the cases 1-3 in Table~\ref{t6}. If $W \cap V = \varnothing$ then the equation of the curve is $y^n=G(x)$ where
\begin{equation*}
\begin{split}
G(x)=\prod_{i=1}^\delta G_{\lambda_i}(x)
\end{split}
\end{equation*}
and $\d$ is as case 1 in Table~\ref{t2}. Let $a_1,...,a_\d$ denote the symmetric polynomials in $\lambda_1,...,\lambda_\d$. Further we can take $\lambda_1...\lambda_\d=1$. Hence the equation of the curve is
\begin{equation*}
\begin{split}
y^n= x^{m\d}+a_1x^{m(\d-1)}+...+a_\delta x^{m}+1
\end{split}
\end{equation*}
If $V \cap W = \phi^{-1}(q_1)$(i.e. case 2 in Table~\ref{t6}) then we know that the equation is $y^n=\varphi(x).G(x)$. Hence the equation is
\begin{equation*}
\begin{split}
y^n= x^{m\d}+a_1x^{m(\d-1)}+...+a_\delta x^{m}+1
\end{split}
\end{equation*}
where $\d$ is as case 2 in Table~\ref{t2}.
If $V \cap W = \phi^{-1}(q_1)\cup \phi^{-1}(q_2) $ (i.e. case 3 in Table~\ref{t6}) then the equation is $y^n=\varphi(x).\chi(x).G(x)$. Hence
\begin{equation*}
\begin{split}
y^n= x(x^{m\d}+a_1x^{m(\d-1)}+...+a_\delta x^{m}+1)
\end{split}
\end{equation*}
where $\d$ is as case 3 in Table~\ref{t2}.

\subsection{$\G \cong D_{2m}$}  Then, $\phi:\bP \to \bP$ has signature $(2,2,m)$. The branch points of $\phi(x)$ are $\infty$ and $\pm 2$. Let $q_1=\infty$, $q_2=2$ and $q_3=-2$. By Lemma ~\ref{l1}, we know that \[\phi(x)=x^m+\frac{1}{x^m}.\] Since $\phi(x)-2=\frac{(x^m-1)^2}{x^m}$ and $\phi(x)+2=\frac{(x^m+1)^2}{x^m}$, $\varphi(x)=x$, $\chi(x)=x^m-1$ and $\psi=x^m+1$. In this case we have $G(x)$ as below.
\begin{equation*}
\begin{split}
G(x)=\prod_{i=1}^\delta (x^{2m}-\lambda_ix^m+1)
\end{split}
\end{equation*}
where $\lambda_i \in $ S  $\backslash$  $\{0,\pm 2\}$ and $\d$ is as corresponding case in Table ~\ref{t2}. Then each family is parameterized as cases 4-9 in Table ~\ref{t6}.

\subsection{$\G \cong A_4$} Then, $\phi:\bP \to \bP$ has signature $(2,3,3)$. We choose branch points $q_1=\infty$, $q_2=6i\sqrt{3}$, and $q_3=-6i\sqrt{3}$, where $i^2=-1$. We know that \[\phi(x)=\frac{x^{12}-33x^8-33x^4+1}{x^2(x^4-1)^2}.\] Thus the points in the fiber of $q_1,q_2,q_3$ are the roots of the polynomials:
\begin{align*}
\begin{split}
\varphi(x)&=x(x^4-1)\\
\chi(x)&=x^4-2i\sqrt{3}x^2+1\\
\psi(x)&=x^4+2i\sqrt{3}x^2+1
\end{split}
\end{align*}
Let $\lambda_i\in $ S  $\backslash$  $\{\infty,\pm 6i\sqrt{3}\}$ then points of $\phi^{-1}(\lambda_i)$ are roots of the polynomial
\begin{equation*}
\begin{split}
G_{\lambda_i}(x)=x^{12}-\lambda_ix^{10}-33x^8+2\lambda_ix^6-33x^4-\lambda_ix^2+1
\end{split}
\end{equation*}
There are $\d$ points in S  $\backslash$  $\{\infty,\pm 6i\sqrt{3}\}$. Hence, we have
\begin{equation*}
\begin{split}
G(x)=\prod_{i=1}^\delta(x^{12}-\lambda_ix^{10}-33x^8+2\lambda_ix^6-33x^4-\lambda_ix^2+1)
\end{split}
\end{equation*}
Then, each family is parameterized as cases 10-15 in Table ~\ref{t6}, where $\d$ is as corresponding case in Table ~\ref{t2}.

\subsection{$\G \cong S_4$} Then, $\phi:\bP \to \bP$ has signature $(2,3,4)$. The branch points of $\phi(x)$ are $\{0,1,\infty\}$. Let $q_1=1$, $q_2=0$ and $q_3=\infty$. Then
\begin{align*}
\begin{split}
\varphi(x)&=x^{12}-33x^8-33x^4+1\\
\chi(x)&=x^8+14x^4+1\\
\psi(x)&=x(x^4-1)
\end{split}
\end{align*}
For $\lambda_i\in $ S  $\backslash$  $\{0,1,\infty\}$, the points in $\phi^{-1}(\lambda_i)$ are roots of the polynomial
\begin{equation*}
\begin{split}
G_{\lambda_i}(x) = & x^{24}+\lambda_ix^{20}+(759-4\lambda_i)x^{16}+2(3\lambda_i+1228)x^{12}\\
 & +(759-4\lambda_i)x^8+\lambda_ix^4+1
\end{split}
\end{equation*}
There are $\d$ points in $S \backslash \{0,1,\infty\}$, where $\d$ is given as in Table ~\ref{t2}. We denote
\begin{equation*}
\begin{split}
M(x):=\prod_{i=1}^\d G_{\lambda_i}(x)
\end{split}
\end{equation*}
Then, each family is parameterized as cases 16-23, where $R(x),S(x),T(x)$ are $\varphi(x),\chi(x),\psi(x)$ respectively.

\subsection{$\G \cong A_5$} The branch points of $\phi:\bP \to \bP$ are $0$, 1728 and $\infty$. Let $q_1=0$, $q_2=\infty$ and $q_3=1728$. At the place $q_3=1728$ the function has the following ramification

\[
\phi(x)-1728=-\frac{(x^{30}+522x^{25}-10005x^{20}-10005x^{10}-522x^5+1)^2}{x^5(x^{10}+11x^5-1)^5}
\]
Then,
\begin{align*}
\begin{split}
\varphi(x)&=x^{20}-228x^{15}+494x^{10}+228x^5+1\\
\chi(x)&=x(x^{10}+11x^5-1)\\
\psi(x)&=x^{30}+522x^{25}-10005x^{20}-10005x^{10}-522x^5+1
\end{split}
\end{align*}
For each $\lambda_i \in$ S $\backslash$  $\{0,1728,\infty\}$ the places in $\phi^{-1}(\lambda_i)$ are the roots of the following polynomial
\begin{equation*}
\begin{split}
G_{\lambda_i}(x)= & -x^{60}+(684-\lambda_i)x^{55}-(55\lambda_i+157434)x^{50}-(1205\lambda_i-12527460)x^{45}\\
                &-(13090\lambda_i+77460495)x^{40}+(130689144-69585\lambda_i)x^{35}\\
                & +(33211924-134761\lambda_i)x^{30}+(69585\lambda_i-130689144)x^{25}\\
                & -(13090\lambda_i+77460495)x^{20}-(12527460-1205\lambda_i)x^{15}\\
                &  -(157434+55\lambda_i)x^{10} +(\lambda_i-684)x^5-1
\end{split}
\end{equation*}
Then,
\begin{equation*}
\begin{split}
\Lambda(x)=\prod_{i=1}^\d G_{\lambda_i}(x)
\end{split}
\end{equation*}
Then equations of the curves are as in cases 24-31 in Table ~\ref{t6}, where 
$Q(x)= \psi (x)$. 

\subsection{$\G \cong U$} The branch point of the curve $\p$ is $\{\infty\}$. Let $q_1=\infty$. Then $\varphi(x)=1$.
For each $\lambda_i \in S\backslash \{\infty\}$ we have
\begin{equation*}
\begin{split}
G_{\lambda_i}(x)=\displaystyle{\prod_{a \in H_t}} (x+a)-\lambda_i
\end{split}
\end{equation*}
There are $\d$ points in $S\backslash \{\infty\}$. Where $\d$ is as in Table ~\ref{t2}. We
denote
\begin{equation*}
\begin{split}
B(x)=\prod_{i=1}^\delta G_{\lambda_i}(x)
\end{split}
\end{equation*}
Then, each family is parameterized as cases 32-33.

\subsection{$\G \cong K_m$} The branch points of the curve $\p$ are $\{0, \infty \}$. Let $q_1=0$, $q_2=\infty$. Then the polynomial over the branch point is
\begin{equation*}
\begin{split}
\varphi(x)&=x\prod_{j=1}^{\frac{p^t-1}{m}}(x^m-b_j)\\
\chi(x)&=1
\end{split}
\end{equation*}
For $\lambda_i \in S\backslash \{0, \infty \}$ we have
\begin{equation*}
\begin{split}
G_{\lambda_i}(x)=((x\prod_{j=1}^{\frac{p^t-1}{m}}(x^m-b_j))^m-\lambda_i)
\end{split}
\end{equation*}
There are $\d$ points in $S\backslash \{0, \infty \}$. Where $\d$ is as in Table ~\ref{t2}. We denote
\begin{equation*}
\begin{split}
\Theta(x)=\prod_{i=1}^\delta G_{\lambda_i}(x)
\end{split}
\end{equation*}
Then, each family is parameterized as cases 34-37.

\subsection{$\G \cong PSL_2(q)$} The branch points of $\p(x)$ are $\{0, \infty \}$. Let $q_1=0$, $q_2=\infty$. Then
\begin{equation*}
\begin{split}
\varphi(x)&=(x^q-x)^{q-1}+1\\
\chi(x)&=x^q-x
\end{split}
\end{equation*}
For $\lambda_i \in S\backslash \{0, \infty \}$, points in $\p^{-1}(\lambda_i)$ are roots of the polynomials,
\begin{equation*}
\begin{split}
G_{\lambda_i}(x)=(((x^q-x)^{q-1}+1)^{\frac{q+1}{2}}-\lambda_i(x^q-x)^{\frac{q(q-1)}{2}})
\end{split}
\end{equation*}
There are $\d$ points in $S\backslash \{0, \infty \}$. Where $\d$ is as in Table ~\ref{t2}. We denote
\begin{equation*}
\begin{split}
\Delta(x)=\prod_{i=1}^{\d}(((x^q-x)^{q-1}+1)^{\frac{q+1}{2}}-\lambda_i(x^q-x)^{\frac{q(q-1)}{2}})
\end{split}
\end{equation*}
Then, each family is parameterized as cases 38-41.

\subsection{$\G \cong PGL_2(q)$} The branch points of $\p(x)$ are $\{0, \infty \}$. Let $q_1=0$, $q_2=\infty$. Then
\begin{equation*}
\begin{split}
\varphi(x)&=(x^q-x)^{q-1}+1\\
\chi(x)&=x^q-x
\end{split}
\end{equation*}
For $\lambda_i \in S\backslash \{0, \infty \}$, points in $\p^{-1}(\lambda_i)$ are roots of the polynomials,
\begin{equation*}
\begin{split}
G_{\lambda_i}(x)=(((x^q-x)^{q-1}+1)^{q+1}-\lambda_i(x^q-x)^{q(q-1)})
\end{split}
\end{equation*}
Then we let,
\begin{equation*}
\begin{split}
\Omega(x)=\prod_{i=1}^{\d}(((x^q-x)^{q-1}+1)^{q+1}-\lambda_i(x^q-x)^{q(q-1)})
\end{split}
\end{equation*}
where $\d$ is gives as Table ~\ref{t2}. Then, each family is parameterized as cases 42-45.
This completes the proof.

\qed


\begin{remark} \label{r2}
By Remark~\ref{r1}, we know that $A_5$ has different ramification when $p=3$. In this case $\phi:\bP \to \bP$ has signature $(6,5)$. The branch points of $\phi(x)$ are $\infty$ and $0$. Let $q_1=\infty$ and $q_2=0$. By Lemma ~\ref{l1}, we know that 
\[\phi(x)=\frac{(x^{10}-1)^6}{(x(x^{10}+2ix^5+1))^5}.\]
Then,

\begin{align*}
\begin{split}
\varphi(x)&=x(x^{10}+2ix^5+1)\\
\chi(x)&=x^{10}-1\\
\end{split}
\end{align*}
For $\lambda_j\in $ S  $\backslash$  $\{0,\infty\}$, the points in $\phi^{-1}(\lambda_j)$ are roots of the polynomial

\begin{equation*}
\begin{split}
G_{\lambda_j}(x)&=x^{60}+\lambda_jx^{55}-(6+10i\lambda_j)x^{50}+35\lambda_jx^{45}+(15+40i\lambda_j)x^{40}+30\lambda_jx^{35}\\
                &\quad-(20-68i\lambda_j)x^{30}+30\lambda_jx^{25}+(15+40i\lambda_j)x^{20}+35\lambda_jx^{15}\\
                &\quad-(6+10i\lambda_j)x^{10}-\lambda_jx^{5}+1
\end{split}
\end{equation*}
There are $\d$ points in S  $\backslash$  $\{0,\infty\}$, where $\d$ is given as in Table ~\ref{t3}. We denote
\begin{equation*}
\begin{split}
P(x):=\prod_{j=1}^\d G_{\lambda_j}(x)
\end{split}
\end{equation*}
Then, each family is parameterized as in Table ~\ref{t7}.
\end{remark}

\medskip

\begin{table}[hb]
\begin{tabular}{|c|c|} 
\hline \hline
$Case$ & $y^n=$ \\
\hline \hline
$a$ &  $P(x)$  \\
$b$ &  $x(x^{10}+2ix^5+1)\cdot P(x)$  \\
$c$ &  $(x^{10}-1)\cdot P(x)$ \\
$d$ &  $x(x^{10}+2ix^5+1)(x^{10}-1)\cdot P(x)$  \\
\hline \hline
\end{tabular} 
\medskip
\caption{Equation of curve when $\G \cong A_5$, $p=3$}\label{t7}
\end{table}

\begin{lemma}
Let $\X_g$  be a cyclic curve defined over an algebraically closed  field $k$ of characteristic $p=3$ such that $\G$ for $\X_g$ is isomorphic to $A_5$. Then, the equation of $\X_g$ is as in one of the cases in Table~\ref{t7}.  

\end{lemma}

We summarize all the cases in the following Theorem.

\begin{thm}
Let $\X_g$ be e genus $g \geq 2$ algebraic curve defined over an algebraically closed field $k$, $G$ its automorphism group over $k$, and $H$  cyclic normal subgroup of   $G$ of  order $n$  such that $g (X_g^H ) =0$. Then, the equation of $\X_g$ can be written as in one of the following cases:

\tiny

\begin{table}[hb]
\begin{center}
\begin{tabular}{|c|c|c|}
\hline
$\#$ &$\G$& $y^n=f(x)$ \\
\hline
1 &     & $x^{m\d}+a_1x^{m(\d-1)}+...+a_\d x^{m}+1$\\
2 &$C_m$& $x^{m\d}+a_1x^{m(\d-1)}+...+a_\d x^{m}+1$\\
3 &     & $x(x^{m\d}+a_1x^{m(\d-1)}+...+a_\d x^{m}+1)$\\
\hline
4 &        & $F(x):= \prod_{i=1}^\d(x^{2m}+\lambda_ix^m+1)$\\
5 &        & $(x^m-1)\cdot F(x)$\\
6 &        & $x\cdot F(x)$\\
7 &$D_{2m}$& $(x^{2m}-1)\cdot F(x)$\\
8 &        & $x(x^m-1)\cdot F(x)$\\
9 &        & $x(x^{2m}-1)\cdot F(x)$\\
\hline
10 &     & $G(x):= \prod_{i=1}^\delta(x^{12}-\lambda_ix^{10}-33x^8+2\lambda_ix^6-33x^4-\lambda_ix^2+1)$\\
11 &     & $(x^4+2i\sqrt{3}x^2+1)\cdot G(x)$\\
12 &$A_4$& $(x^8+14x^4+1)\cdot G(x)$\\
13 &     & $x(x^4-1)\cdot G(x)$\\
14 &     & $x(x^4-1)(x^4+2i\sqrt{3}x^2+1)\cdot G(x)$\\
15 &     & $x(x^4-1)(x^8+14x^4+1)\cdot G(x)$\\
\hline
16 &     & $M(x)$\\
17 &     & $  \left( x^8+14x^4+1 \right)  \cdot M(x)$\\
18 &     & $x(x^4-1) \cdot M(x)$\\
19 &     & $\left( x^8+14x^4+1 \right)  \cdot x(x^4-1) \cdot M(x)$\\
20 &$S_4$& $\left( x^{12}-33x^8-33x^4+1  \right)\cdot M(x)$\\
21 &     & $\left( x^{12}-33x^8-33x^4+1  \right)  \cdot \left( x^8+14x^4+1 \right)  \cdot M(x)$\\
22 &     & $\left( x^{12}-33x^8-33x^4+1  \right) \cdot x(x^4-1) \cdot M(x)$\\
23 &     & $\left( x^{12}-33x^8-33x^4+1  \right) \cdot \left( x^8+14x^4+1 \right)  \cdot x(x^4-1)  M(x)$\\
\hline
24 &     & $\Lambda(x)$\\
25 &     & $(x(x^{10}+11x^5-1))\cdot \Lambda(x)$\\
26 &     & $(x^{20}-228x^{15}+494x^{10}+228x^5+1)(x(x^{10}+11x^5-1))\cdot \Lambda(x)$\\
27 &     & $(x^{20}-228x^{15}+494x^{10}+228x^5+1)\cdot \Lambda(x)$\\
28 &$A_5$& $Q (x) \cdot \Lambda(x)$\\
29 &     & $x(x^{10}+11x^5-1).\psi(x)\cdot \Lambda(x)$\\
30 &     & $(x^{20}-228x^{15}+494x^{10}+228x^5+1)\cdot \psi(x)\cdot\Lambda(x)$\\
31 &     & $(x^{20}-228x^{15}+494x^{10}+228x^5+1)(x(x^{10}+11x^5-1))\cdot\psi(x)\cdot\Lambda(x)$\\
\hline
32 &$U$& $B(x)$\\
33 &   & $B(x)$\\
\hline
34 &     & $\Theta(x)$\\
35 &$K_m$& $x\prod_{j=1}^{\frac{p^t-1}{m}}\left(x^m-b_j\right)\cdot\Theta(x)$\\
36 &     & $\Theta(x)$\\
37 &     & $x\prod_{j=1}^{\frac{p^t-1}{m}}\left(x^m-b_j\right)\cdot\Theta(x)$\\
\hline
38 &          & $\Delta(x)$\\
39 &$PSL_2(q)$& $((x^q-x)^{q-1}+1)\cdot\Delta(x)$\\
40 &          & $(x^q-x)\cdot\Delta(x)$\\
41 &          & $(x^q-x)((x^q-x)^{q-1}+1)\cdot\Delta(x)$\\
\hline
42 &          & $\Omega(x)$\\
43 &$PGL_2(q)$& $((x^q-x)^{q-1}+1)\cdot\Omega(x)$\\
44 &          & $(x^q-x)\cdot\Omega(x)$\\
45 &          & $(x^q-x)((x^q-x)^{q-1}+1)\cdot\Omega(x)$\\
\hline
\end{tabular}
\medskip
\caption{The equations of the curves related to the cases in Table ~\ref{t2}}
\end{center}
\end{table}

\normalsize

where   $\d$ is given as in Table~\ref{t2} and $M, \Lambda, Q, B, \Delta,$ and $\Omega$ are as follows: 

\begin{equation*}
\begin{split}
M = & \prod_{i=1}^\d \left(  x^{24}+\lambda_ix^{20}+(759-4\lambda_i)x^{16}+2(3\lambda_i+1228)x^{12} \right. \\
   & + \left. (759-4\lambda_i)x^8+\lambda_ix^4+1 \right) \\
\Lambda= & \prod_{i=1}^\d \left(-x^{60}+(684-\lambda_i)x^{55}-(55\lambda_i+157434)x^{50}-(1205\lambda_i-12527460)x^{45}\right.\\
                &-(13090\lambda_i+77460495)x^{40}+(130689144-69585\lambda_i)x^{35}\\
                & +(33211924-134761\lambda_i)x^{30}+(69585\lambda_i-130689144)x^{25}\\
                & -(13090\lambda_i+77460495)x^{20}-(12527460-1205\lambda_i)x^{15}\\
                &  \left. -(157434+55\lambda_i)x^{10} +(\lambda_i-684)x^5-1 \right) \\
Q  = & x^{30}+522x^{25}-10005x^{20}-10005x^{10}-522x^5+1\\
B = &  \prod_{i=1}^\delta \displaystyle{\prod_{a \in H_t}}   \left( (x+a)-\lambda_i \right) \\
\Delta = & \prod_{i=1}^{\d}(((x^q-x)^{q-1}+1)^{\frac{q+1}{2}}-\lambda_i(x^q-x)^{\frac{q(q-1)}{2}}) \\
\Omega = & \prod_{i=1}^{\d}(((x^q-x)^{q-1}+1)^{q+1}-\lambda_i(x^q-x)^{q(q-1)}) \\
\end{split}
\end{equation*}

\end{thm}

\end{document}